\documentclass[12pt]{amsart}   

\usepackage{graphicx}
\usepackage{tikz}
\usepackage{amssymb}
\usetikzlibrary{calc}
\usepackage{url}
\usepackage{caption}
\usepackage{tcolorbox}
\usepackage{listings}

\newtheorem{thm}{Theorem}[section]

\theoremstyle{definition}

\theoremstyle{remark}

\begin{document}

\title{A Graph Theoretic Formula \\for the Number of Primes $\pi(n)$
}

\author{R. Jacobs, C. E. Larson$^*$}


\address{Department of Mathematics and Applied Mathematics\\Virginia Commonwealth University\\Richmond, VA 23284, USA }

\thanks{(*) Research supported by the Simons Foundation Mathematics and Physical Sciences--Collaboration Grants for Mathematicians Award (426267)}


\date{}

\maketitle

\begin{abstract}
Let PR$[n]$ be the graph whose vertices are $2,3,\ldots,n$ with vertex $v$ adjacent to vertex $w$ if and only if $\gcd(v,w)>1$. It is shown that $\pi(n)$, the the number of primes no more than $n$, equals the Lov\'{a}sz number of this graph. This result suggests new avenues for graph-theoretic investigations of number-theoretic problems.
\end{abstract}

In \textit{Written on the Wall} (or \textit{WoW}), Fajtlowicz's notes on conjectures of his program \textsc{graffiti} \cite{Fajt88a,Fajt95}, Fajtlowicz defined the graphs RP$[S]$ and PR$[S]$ whose vertices are a set $S$ of integers \cite{Fajt}. For RP$[S]$ two distinct vertices are adjacent if and only if they are relatively prime; while in PR$[S]$ two distinct vertices are adjacent if and only if they have a non-trivial common factor (and are thus the complements of the $RP[S]$ graphs). Let PR$[n]$ be the graph where $S=\{2,3,\ldots,n\}$. \textit{WoW} is indexed by conjecture numbers, often with useful commentary of its author and correspondents. These graphs are defined in WoW \#434 (1988).  Study of these graphs may yield new insights into number theoretic questions. Among other things WoW records  Staton's proof of the interesting fact that \textit{every graph is an induced subgraph of PR$[n]$, for some positive integer $n$} (WoW \#446).

An independent set in a graph is a set of vertices which are pair-wise nonadjacent. The independence number $\alpha=\alpha(G)$ of a graph $G$ is the cardinality of a maximum independent set. 
Fajtlowicz observed, and it is easy to see, that the primes are a maximum independent set in PR$[n]$ and thus the independence number of PR$[n]$ is the number of primes up to $n$, denoted $\pi(n)$. So $\alpha(\text{PR}[n])=\pi(n)$. The Prime Number Theorem gives an asymptotic formula for $\pi(n)$ and the Riemann Hypothesis is a equivalent to a conjecture about the error term in a formula for $\pi(n)$. 
Thus study of the independence number of the PR$[n]$ graphs may yield new insights into $\pi(n)$. \textsc{graffiti} is well-known for its conjectures for the independence number of a graph---many of which were proved. 
While following up on Fajtlowicz's idea and also investigating the class of graphs where $\alpha$ equals Lov\'{a}sz's theta function $\vartheta$ we discovered:

\begin{thm}
$\pi(n)=\vartheta(\text{PR}[n])$
\end{thm}

This may be of interest for a few reasons: Lov\'{a}sz's theta function is widely studied, has a number of equivalent definitions \cite{Knut94}, is efficiently computable, and there is an efficient algorithm for recognizing graphs where the independence number $\alpha$ equals Lov\'{a}sz's theta function $\vartheta$. 

An \textit{orthonormal representation} of a graph $G$ is an assignment of a unit vector $\hat{x}_v$ to each vertex $v$ in the vertex set $V(G)$ having the property that vectors assigned to non-adjacent vertices are orthogonal. The Lov\'{a}sz's theta function $\vartheta=\vartheta(G)$ of a graph $G$ is defined to be a minimum over all feasible orthonormal representations $\{ \hat{x}_v: v\in V(G)\}$ (or simply $\{\hat{x}_v\}$) and all unit vectors $\hat{c}$:
\[
\vartheta = \min_{\hat{c},\{\hat{x}_v\}} \max_{v\in V(G)} \frac{1}{(\hat{c}\cdot \hat{x}_v)^2}.
\]

\begin{proof}
Let $n$ be a positive integer. Since Lov\'{a}sz proved that for any graph $G$, $\vartheta(G)\geq \alpha(G)$ \cite{Lova79} and $\alpha(PR[n])=\pi(n)$, it is enough to show that there is in fact a feasible orthonormal representation of PR$[n]$ and unit vector $\hat{c}$ that realizes this lower bound; that is, such that:
\[
\pi(n) = \min_{\hat{c},\{\hat{x}_v\}} \max_{v\in V(G)} \frac{1}{(\hat{c}\cdot \hat{x}_v)^2}.
\]
Suppose there are exactly $k$ primes no more than $n$: $p_1,p_2,\ldots,p_k$. 
For each vertex $v$ with $l_v$ distinct prime factors we define $\hat{x}_v$ to be a $k$-component vector where the $i^{th}$ component is $0$ unless $p_i$ is a factor of $v$ in which case the $i^{th}$ component is $\frac{1}{\sqrt{l_v}}$. 
It is easy to see that if $v$ and $w$ are relatively prime, and thus non-adjacent, then $\hat{x}_v \cdot \hat{x}_w=0$. 
Then let $\hat{c}$ be the vector whose components are all $\frac{1}{\sqrt{k}}$. 
Then $\max_{v\in V(G)} \frac{1}{(\hat{c}\cdot \hat{x}_v)^2}=\max_{v\in V(G)} \frac{1}{(\l_v\cdot \frac{1}{\sqrt{k}}\cdot \frac{1}{\sqrt{l_v}})^2}=\max_{v\in V(G)}\frac{k}{l_v}$ occurs for a vertex $w$ when $l_w=1$, for instance when $w$ is prime, and in this case equals $k$, the number of primes no more than $n$.
\end{proof}

Since perfect graphs have the property that $\alpha=\vartheta$ it might be thought that the PR$[n]$ graphs are perfect. They  are for $n<35$. $PR[35]$ has an odd hole: $5\cdot 7, 7\cdot 3, 3\cdot 11, 11\cdot 2, 2\cdot 5$.

\bibliographystyle{plain}
\bibliography{../../larson.bib}
%

\end{document}